\documentclass{amsart}
\usepackage{amssymb, amsmath, color}
\newtheorem{theorem}{Theorem}

\newtheorem{remark}{Remark}

\theoremstyle{definition}

\newtheorem{corollary}[theorem]{Corollary}

\newcommand{\sech}{\mathrm{sech}}

\begin{document}

\title[Sharp lower bound for the scalar curvature of steady solitons]
{A sharp lower bound for the scalar curvature of certain steady gradient Ricci solitons}
\author{Manuel Fern\'{a}ndez-L\'{o}pez, Eduardo Garc\'{\i}a-R\'{\i}o}
\address{Faculty of Mathematics, University of Santiago de Compostela,
15782 Santiago de Compostela, Spain} \email{manufl@edu.xunta.es,
eduardo.garcia.rio@usc.es}
\thanks{Supported by projects MTM2006-01432 and PGIDIT06PXIB207054PR (Spain)}



\begin{abstract}
We show that the scalar curvature of a steady gradient Ricci soliton satisfying that the ratio between the square norm of the Ricci tensor and the square of the scalar curvature is bounded by one half, is boundend from below by the hyperbolic secant of one half the distance function from a fixed point. 
\end{abstract}

\maketitle

\section{Introduction}

A \emph{Ricci soliton} is a Riemannian manifold $(M,g)$ that admits a smooth vector field $X$ on $M$ such that
\begin{equation}\label{soliton}
\frac12\mathcal{L}_Xg+ Rc=\lambda g,
\end{equation}
where $\mathcal{L}_X$ is the Lie derivative in the direction of
the vector field $X$, $Rc$ denotes the Ricci tensor and $\lambda$ is a constant. When
the vector field $X$ can be replaced by the gradient of some smooth
function $f$ on $M,$ called the \emph{potential function},
$(M,g)$ is said to be a \emph{gradient Ricci soliton}. In such a case the
equation (\ref{soliton}) becomes
\begin{equation}\label{gradientsoliton}
Rc+H_f=\lambda g,
\end{equation}
where $H_f$ denotes the Hessian of the function $f$. A Ricci soliton (\ref{soliton}) is said to be \emph{shrinking}, \emph{steady} or \emph{expanding} according to $\lambda >0$, $\lambda=0$ or $\lambda <0$. \cite{C} is a very intesting paper for recent information on this topic.

For steady gradient Ricci solitons it is well-known (\cite{C}, for example) that $R+|\nabla f|^2=C,$ where $C$ is a positive constant, unless the steady soliton is Ricci flat. We scale the metric to have the constant equal to $1.$ 

Very recently it was given in \cite{CLY} a lower bound for the scalar curvature of a steady gradient Ricci solitons in terms of the dimension of the manifold and the potential function, under additional assumptions. However there is no a good knwoledge of the behaviour of the potential function of a steady gradient soliton, and thus the bound cannot be expresed in an explicit way in terms of the distance function. In this paper we show the following

\begin{theorem}\label{th}
Let $(M,g)$ be a complete gradient steady Ricci soliton satifying $|Rc|^2\leq \frac{R^2}{2}$ and normalized as before. Then $$R(x)\geq k\sech^2 \frac{r(x)}{2},$$ where $r(x)$ is the distance from a fixed point $O\in M$ and $k\leq 1$ is a constant that only depends on $O$ and $R(O).$
\end{theorem}

\begin{remark}
Note that on the scaled Hamilton's cigar soliton $\left(\mathbb{R}^2, \frac{4(dx^2+dy^2)}{1+x^2+y^2}\right)$ one has $R+|\nabla f|^2=1$ and $R(x)=\sech^2 \frac{r(x)}{2},$ where the distance $r(x)$ is measured from the only point where the scalar curvature attains its maximum. This shows that our lower bound is sharp in dimension two. In higher dimensions, if we consider the product of Hamilton's cigar soliton and any complete Ricci flat manifold, we also have that our bound is sharp. Indeed, note that we actually have equality when moving in the direction of the cigar, where the distance $r(x)$ is measured again from the only point where the scalar curvature attains its maximum.
\end{remark}

It is also possible to give a lower bound for the scalar curvature assuming that the Ricci tensor is nonnegative. In such a case we do not know if such a bound is sharp, because we do not know of any example where the equality is achieved.

\begin{corollary}\label{co}
Let $(M,g)$ be a complete gradient steady Ricci soliton with nonnegative Ricci curvature and normalized as before. Then $$R(x)\geq k\sech^2 r(x),$$ where $r(x)$ is the distance from a fixed point $O\in M$ and $k\leq 1$ is a constant that only depends on $O$ and $R(O).$
\end{corollary}

\section{Proofs of the results}

\noindent{\bf Proof of Theorem \ref{th}.-}

From Kato's inequality, at every point where $\nabla  f$ does not vanish, we have 

$$|H_f|^2=|\nabla \nabla f|^2\geq |\nabla |\nabla f||^2 = |\nabla \sqrt{1-R}|^2=\frac{|\nabla R|^2}{4|\nabla f|^2,}$$ or, equivalently,

$$|\nabla R|^2\leq 4 |H_f|^2|\nabla f|^2.$$ By assumption we have that $|H_f|^2=|Rc|^2\leq \frac{R^2}{2}.$ Thus 

$$\frac{|\nabla R|}{R\sqrt{1 -R}}\leq 1.$$ Now, let $O$ a fixed point on $M$ and let $\gamma:[0,t]\rightarrow M$ be a minimizing geodesic with $\gamma(0)=O.$ Integrating the function $\frac{-(R\circ \gamma)'}{R\sqrt{1 -R}}$ along $\gamma(s)$ we get $$\left[ \ln \frac{1+\sqrt{1-R}}{1-\sqrt{1-R}}\right]_0^t= -\int_0^t  \frac{(R\circ \gamma)'}{R\sqrt{1 -R}}ds \leq \int_0^t  \frac{|\nabla R|^2}{R\sqrt{1 -R}} ds \leq t.$$ Writing $c=\frac{1+\sqrt{1-R(O)}}{1-\sqrt{1-R(O)}}$ we get that $$1+\sqrt{1-R(\gamma(t))}\leq ce^t (1-\sqrt{1-R(\gamma(t))}).$$ Then it is a straightforward computation to get that $$R(\gamma(t))\geq \frac{4c}{c^2e^t+2c+e^{-t}}. $$ Now, since $c\geq 1,$ we have that $$R(\gamma(t))\geq \frac{4c}{c^2e^t+2c+e^{-t}}\geq \frac{4c}{c^2e^t+2c^2+c^2e^{-t}} =\frac{1}{c}\sech^2 \frac{t}{2}.$$ Since the geodesic $\gamma$ and $t$ are arbitrary we have finished the proof.
$\hfill{q.e.d}$

\begin{remark}
Note that the scalar curvature $R$ may take the value $1$ along the geodesic $\gamma([0,t]).$ Since gradient Ricci solitons are analytic manifolds we have two possibilities. If the set of points where $R$ takes the value $1$ on $\gamma([0,t])$ has a point of accumulation then $R$ must be constant and the proof works. In other case, we only have a finite set of points $\gamma([0,t])$ where $R$ takes the value $1.$ Then we have to deal, eventually, with a finite number of improper integrals, and the proof also works. 
\end{remark}

\noindent{\bf Proof of Corollary \ref{co}.-}
Proceeding as in the proof of the theorem and using the inequality $|H_f|^2= |Rc|^2\leq R^2$ we get $$\frac{|\nabla R|}{R\sqrt{1 -R}}\leq 2.$$ Then the result is obtained following the same steps as in the theorem. $\hfill{q.e.d}$

\bibliographystyle{amsplain}

\end{document}